\documentclass[10pt]{amsart}
\usepackage{amssymb}
\usepackage{amsmath}
\usepackage{amsthm}
\usepackage[latin1]{inputenc}

\usepackage[english]{babel}
\usepackage{color, bbm}

\numberwithin{equation}{section}

\newtheorem{thm}{Theorem}[section]
\newtheorem{lem}[thm]{Lemma}

\newtheorem{Rem}[thm]{Remark}

\def \ds{\displaystyle}
\def \e{\epsilon}

\newcommand{\R}{\mathbb{R}}
\newcommand{\B}{\mathbb{B}}
\def \p{\partial}

\def \M{{\mathcal M}}

\begin{document}
	\title{On sharp anisotropic Hardy inequalities}

\author{Xia Huang}
	\address{Xia Huang \newline\indent School of Mathematical Sciences, Key Laboratory of MEA (Ministry of Education) \& Shanghai Key Laboratory of PMMP, East China Normal University, \newline\indent
		Shanghai, 200241, People's Republic of China}
	\email{ xhuang@cpde.ecnu.edu.cn}
\author{Dong Ye}
	\address{Dong Ye\newline\indent School of Mathematical Sciences, Key Laboratory of MEA (Ministry of Education) \&  Shanghai Key Laboratory of PMMP, East China Normal University, \newline\indent
		Shanghai, 200241, People's Republic of China}
	\email{ dye@math.ecnu.edu.cn}

\begin{abstract}
Recently, Yanyan Li and Xukai Yan showed in \cite{LY,LY2} the following interesting Hardy inequalities with anisotropic weights: Let $n\geq 2$, $p \geq 1$, $p\alpha > 1-n$, $p(\alpha + \beta)> -n$, then there exists $C > 0$ such that
$$\||x|^{\beta}|x'|^{\alpha+1} \nabla u\|_{L^p(\mathbb{R}^n)} \geq C\||x|^\beta|x'|^\alpha u\|_{L^p(\mathbb{R}^n)}, \quad \forall\; u\in C_c^1(\mathbb{R}^n).$$
Here $x' = (x_1,\ldots, x_{n-1}, 0)$ for $x = (x_i) \in \R^n$. In this note, we will determine the best constant for the above estimate when $p=2$ or $\beta \geq 0$. Moreover, as refinement for very special case of Li-Yan's result in \cite{LY2}, we provide explicit estimate for the anisotropic $L^p$-Caffarelli-Kohn-Nirenberg inequality.
\end{abstract}
\maketitle

\noindent
{\footnotesize {\bf 2020 MSC:} {\sl 26D10, 35A23, 46E35}}

\noindent
{\footnotesize {\bf Keywords:} {\sl Anisotropic Hardy inequality, best constant, anisotropic $L^p$-Caffarelli-Kohn-Nirenberg inequality}}
	
\section{Introduction}

\medskip
Recently, Li-Yan studied in \cite{LY} the asymptotic stability of the $(-1)$-homogeneous axisymmetric stationary solutions to Navier-Stokes equations in $\R^3$. A key estimate was
\begin{equation}\label{0.1}
\int_{\mathbb{R}^3}|\nabla u|^2 dx \geq C\int_{\mathbb{R}^3}\frac{u^2}{|x'||x|}dx, \quad \forall \; u \in C_c^1(\R^3),
\end{equation}
where $x' = (x_1, x_2, 0)$. Li-Yan \cite{LY} showed \eqref{0.1} by proving the following strengthened inequality:
\begin{equation}\label{0.2}
\int_{\mathbb{R}^3} \frac{|x'|}{|x|}|\nabla u|^2 dx \geq C\int_{\mathbb{R}^3}\frac{u^2}{|x'||x|}dx, \quad  \forall \; u \in C_c^1(\R^3).
\end{equation}
It is worthy to remark that \eqref{0.2} improves also the classical Hardy inequality
\begin{equation*}
\int_{\mathbb{R}^3}|\nabla u|^2 dx \geq \frac{1}{4}\int_{\mathbb{R}^3}\frac{u^2}{|x|^2} dx, \quad  \forall \; u \in C_c^1(\R^3).
\end{equation*}
This motivated them to show the following general anisotropic Hardy inequalities, see \cite[Theorem 1.3]{LY} and \cite[pages 6-7]{LY2}.

\smallskip
\noindent
{\bf Theorem A.} {\it Let $n\geq 2$, $x' = (x_1,\ldots, x_{n-1}, 0)$ for $x = (x_i) \in \R^n$. Assume that $p\geq 1$, $p\alpha >1-n$, $p(\alpha+\beta)>-n$, then there exists positive constant $C$ depending on $n, p, \alpha$ and $\beta$, such that
\begin{equation}\label{ly}
\left\||x|^{\beta}|x'|^{\alpha+1} \nabla u\right\|_{L^p(\mathbb{R}^n)} \geq C\left\||x|^\beta|x'|^\alpha u\right\|_{L^p(\mathbb{R}^n)}, \quad \forall\; u\in C_c^1(\mathbb{R}^n).
\end{equation}}

The interesting estimates \eqref{ly} were used by Li-Yan to establish a generalized and improved anisotropic version of Caffarelli-Kohn-Nirenberg's interpolation inequalities, see \cite[Theorem 1.1]{LY2}, which established necessary and sufficient conditions to have
\begin{align}
\label{CKN}
\left\||x|^{\gamma_1}|x'|^\alpha u\right\|_{L^s(\mathbb{R}^n)} \leq C \left\||x|^{\gamma_2}|x'|^{\mu} \nabla u\right\|^a_{L^p(\mathbb{R}^n)}\left\||x|^{\gamma_3}|x'|^\beta u\right\|^{1-a}_{L^q(\mathbb{R}^n)}, \quad \forall\; u \in C_c^1(\R^n).
\end{align}
%In \cite{BC}, Bao-Chen consider the existence and symmetry  breaking region of extremal functions, and their symmetry properties for equality \eqref{CKN} with $p=2$ and $q=1$.

\medskip
Our main purpose here is to study the best constant for the inequalities \eqref{ly}. A first remark is that for any $n\geq 2$, $p \geq 1$, $|x|^\beta|x'|^\alpha \in L^p_{loc}(\R^n)$ if and only if
\begin{align}
\label{cond1.1}
p\alpha > 1-n \quad \mbox{and}\quad p(\alpha + \beta)> -n.\end{align}
This fact can be seen by spherical coordinates. Our first result gives a complete answer for best constant of \eqref{ly} when $p=2$.
\begin{thm}
\label{thm1.1}
Let $n\geq 2$, $p = 2$, and $\alpha, \beta\in\mathbb{R}$ satisfy \eqref{cond1.1}. Denote by $C_{n,\alpha, \beta}$ the sharp constant in \eqref{ly} with $p=2$, i.e., the best constant to claim
%$K:=-4\beta(n+2\alpha+\beta)\geq 0$, Then for any {\color{red} $u\in C_c^1(\mathbb{R}^n)$}, there holds
\begin{equation}
\label{p=2}
\int_{\mathbb{R}^n}|x'|^{2\alpha+2}|x|^{2\beta}|\nabla u|^2 dx \geq C_{n, \alpha, \beta}\int_{\mathbb{R}^n} |x'|^{2\alpha}|x|^{2\beta}u^2 dx, \quad \forall\; u\in C_c^1(\mathbb{R}^n).
\end{equation}
Then we have
\begin{equation}
\label{const1.1} C_{n, \alpha, \beta}   =
\frac{(n-1+2\alpha)^2-\big[\sqrt{\max(K, 1)} - 1\big]^2}{4},
\end{equation}
where $K = -4\beta(n+2\alpha+\beta)= (n+2\alpha)^2 - (n+2\alpha+2\beta)^2$.
\end{thm}

\medskip
For more general $p \geq 1$, we obtain the following partial result where the best constant is determined when $\beta \geq 0$.
\begin{thm}
\label{thmp}
Let $n\geq 2$, $p \geq 1$ and $\alpha, \beta \in \R$ satisfy \eqref{cond1.1}. Denote still by $C_{n,\alpha,\beta}$ the best constant in \eqref{ly}, then
$$C_{n,\alpha,\beta}=\Big(\frac{n-1+ p\alpha}{p}\Big)^p \quad \mbox{for any }\; \beta \geq 0.$$
Moreover, $C_{n,\alpha,\beta}\geq \Big(\frac{n-1+ p\alpha + p\beta}{p}\Big)^p$ if $\beta<0$ and $p(\alpha+\beta)>1-n$.
\end{thm}

Furthermore, we give an alternative proof for the anisotropic $L^p$-Caffarelli-Kohn-Nirenberg inequalities, i.e., a very special case of Li-Yan's general result \eqref{CKN} with $s=p=q > 1$ and $a = \frac{1}{p}$. Let $n\geq 2$, $p\geq 1$, we know that $|x'|^\mu|x|^{\gamma_2}$, $|x'|^\beta|x|^{\gamma_3}$, $|x'|^{\alpha}|x|^{\gamma_1} \in L_{loc}^p(\mathbb{R}^n)$ if and only if
\begin{align}\label{cond1.2}
\min(\alpha, \beta, \mu) >  \frac{1-n}{p}, \quad \min(\alpha+\gamma_1, \mu+\gamma_2, \beta+\gamma_3)>-\frac{n}{p}.
\end{align}
 Notice also that in our special case, the assumptions \cite[(1.10)-(1.13)]{LY2} are equivalent to
\begin{align}
\label{condLY}
\alpha + \gamma_1 = \frac{\mu+ \gamma_2-1}{p} + \frac{(p-1)(\beta + \gamma_3)}{p}, \quad \gamma_1 \le \frac{\gamma_2-1}{p} + \frac{(p-1)\gamma_3}{p}.
\end{align}
\begin{thm}\label{thm1.3}
Let $n\geq 2$, $p>1$. Assume that $\alpha, \beta, \mu, \gamma_1, \gamma_2, \gamma_3$ satisfy \eqref{cond1.2}-\eqref{condLY}, then for any $u \in C_c^1(\R^n)$, there holds
\begin{align}
\label{speCKN}
\left\||x|^{\gamma_2}|x'|^{\mu} \nabla u\right\|_{L^p(\mathbb{R}^n)}\left\||x|^{\gamma_3}|x'|^\beta u\right\|^{p-1}_{L^p(\mathbb{R}^n)} \geq \frac{n+p(\alpha+\gamma_1)}{p}\left\||x|^{\gamma_1}|x'|^\alpha u\right\|^p_{L^p(\mathbb{R}^n)}.
\end{align}
If moreover $\alpha = \beta=\mu$ and $\gamma_3 - \gamma_2 + 1 > 0$, the constant $\frac{n+p(\alpha+\gamma_1)}{p}$ is sharp.
\end{thm}

Our approach departs from an elementary identity: Let $\Omega \subset \R^n$ be an open set, $V \in C^1(\Omega)$ and $f\in C^2(\Omega)$ be positive, then
\begin{align}
\label{E2}
\int_{\Omega} V |\nabla u|^2 dx  = - \int_{\Omega} \frac{{\rm div}(V \nabla f)}{f}u^2 dx + \int_{\Omega} Vf^2\Big|\nabla\Big(\frac{u}{f}\Big) \Big|^2 dx, \quad \forall\; u \in C^1_c(\Omega).
\end{align}
The above equality can be showed by using integration by parts; or by taking $\vec F=-\frac{\nabla f}{f}$ in the more general equality
\begin{align}
\label{E1}
\int_{\Omega} V |\nabla u|^2 dx = \int_{\Omega} \left[{\rm div}(V \vec{F})-V|\vec {F}|^2\right]u^2 dx + \int_{\Omega} V\big|\nabla u + u\vec {F}\big|^2 dx.
\end{align}
These identities suggest to find weighted Hardy or Poincar\'e inequalities by testing suitable positive functions $f \in C_c^1(\Omega)$, and provide a natural way to study Hardy type inequalities. This idea has been used in many situations in the literature, and summarized in \cite{HY}. For example, the Bessel pair with radial potential $(V, W)$ introduced by Ghoussoub-Moradifam \cite{GM} is a special case of \eqref{E2} for radial function $f(x)=f(|x|)$, since
$$f''(r) + \Big(\frac{n-1}{r}+ \frac{V'}{V}\Big)f'(r) + \frac{cW}{V}f(r) = 0 \;\; \mbox{is equivalent to say } {-{\rm div}(V\nabla f)} = cWf.$$

Furthermore, we remark that the last integral in \eqref{E2} is zero if and only if $u/f$ is a constant. It is well known that the optimal Hardy inequality cannot be reached in general, that is, the {\sl best} choice of $f$ does not belong to the proper functional space. However, we can check eventually sharpness of the subsequent weight $W := -\frac{{\rm div}(V \nabla f)}{f}$ by choosing appropriate functions $u$ to {\sl approximate} $f$.

\smallskip
The last term in \eqref{E2} can be interpreted also as a kind of stability, since it measures in some sense the {\sl distance} between $u$ and the eventual linear space generated by the optimal choice $f$.

\medskip
To prove Theorem \ref{thm1.1}, according to $V$, we will apply \eqref{E2} with $f(x) = |x'|^\theta|x|^\lambda$, and try to optimize the subsequent weight $W$ with suitable choice of the parameters $\theta$, $\lambda \in \R$. As $\theta$ or $\lambda$ are allowed to be negative, the corresponding anisotropic Hardy inequalities are firstly proved in $C_c^1(\R^n\backslash\{x' = 0\})$, then extended to $C_c^1(\R^n)$ by density argument. Moreover, we study the sharpness by trying to approximate the {\sl optimal} choice of $f$.

\medskip
%\footnote{\; The case $p=1$ requires more assumption on $f$, so need to be explained a little bit}
An equality similar to \eqref{E2} exists for general $p > 1$, where we replace the last integral by a Picon\'e type term, see \cite[section 10.2]{HY}. Let $(\mathcal{M}, g)$ be a Riemannian manifold, consider $V\in C^1(\M)$ and $\vec F\in C^1(\M, T_g\M)$, then for any $u\in C_c^1(\M)$, there holds
\begin{align}
\label{ELpF}
\int_{\M}[ V|\nabla u|^p +(p-1)V|\vec F|^p|u|^p ]dg= \int_{\M} {\rm div}\left(V|\vec F|^{p-2}\vec F\right)|u|^p dg+\int_\M V \mathcal{R}(\nabla u, u\vec F)dg,
\end{align}
where
\begin{align*}
{\mathcal R}(\vec X, \vec Y) & = (p-1)\|\vec Y\|^p + \|\vec X\|^p + p\|\vec Y\|^{p-2} \langle\vec Y, \vec X \rangle \geq 0.
\end{align*}
In particular, let $\vec F=-\frac{\nabla f}{f}$ and $\M = \Omega$, we get that (see also \cite[Theorem 3.1]{MN}, \cite{PT} for $V \equiv 1$, or \cite{DFLL, DLL, FLL} with $f$ depending on one variable)
\begin{align}
\label{Ep}
\begin{split}
\int_{\Omega} V|\nabla u|^p dx & = - \int_{\Omega} \frac{{\rm div}\left(V|\nabla f|^{p-2}\nabla f\right)}{f^{p-1}}|u|^p dx + \int_{\Omega} V{\mathcal R}\Big(\nabla u, -u\frac{\nabla f}{f}\Big)dx.
\end{split}
\end{align}
Hence we obtain the $L^p$-Hardy inequality (with suitable $V$, $f$ and $u$)
\begin{align}
\label{HLp}
\int_{\Omega} V|\nabla u|^p dx \geq -\int_{\Omega} \frac{{\rm div}\left(V|\nabla f|^{p-2}\nabla f\right)}{f^{p-1}}|u|^p dx.
\end{align}
Here again, we need not any symmetry assumption on $V$ or $f$, again the residual term in \eqref{Ep} is zero if and only if $u/f$ is constant. Therefore we can proceed similarly as for $L^2$ case to handle Theorem \ref{thmp}.
\begin{Rem}
The identity \eqref{HLp} holds also for $p=1$, if $V \in C^1(\Omega)$, $f\in C^{1,1}(\Omega)$ satisfies $|\nabla f|>0$. In that case,
\begin{align*}
\int_{\Omega} V|\nabla u|dx \geq -\int_{\Omega} {\rm div}\left(V\frac{\nabla f}{|\nabla f|}\right)|u| dx, \quad \forall\; u\in C_c^1(\Omega).
\end{align*}
\end{Rem}

\medskip
Moreover, notice that for any $\kappa > 0$, $(\kappa^{-1}V, \kappa^{\frac{1}{p-1}}\vec F)$ does not change the subsequent weight $W = {\rm div}(V|\vec F|^{p-2}\vec F)$ on the right hand side of \eqref{ELpF}. Taking
\begin{align*}
\kappa_0^\frac{p}{p-1} = \frac{\ds \int_{\Omega}V|\nabla u|^p dg}{\ds \int_{\Omega} V|\vec F|^p u^p dg},
\end{align*}
we obtain a special weighted $L^p$-Caffarelli-Kohn-Nirenberg inequality as follows.
\begin{align}
\label{CKNp}
\begin{split}
& \quad \left(\int_{\Omega}V|\nabla u|^p dg\right)^{\frac{1}{p}} \left(\int_\Omega V|\vec F|^p|u|^p dg\right)^{\frac{p-1}{p}}\\ & = \frac{1}{p}\int_{\Omega} {\rm div}(V|\vec F|^{p-2}\vec F)|u|^p dg+\int_\omega \frac{V}{p\kappa_0} \mathcal{R}\Big(\nabla u, u\kappa_0^{\frac{1}{p-1}}\vec F\Big)dg\\
& \geq \frac{1}{p}\int_{\Omega} {\rm div}(V|\vec F|^{p-2}\vec F)|u|^p dg.
\end{split}
\end{align}
We will choose suitable $\vec F$ to prove Theorem \ref{thm1.3}.

\begin{Rem}
Very recently, Li-Yan's anisotropic Hardy inequalities are generalized by Musina-Nazarov \cite{MN}, see section 5. For Li-Yan's inequality \eqref{CKN} with $p=2$ and $a=1$, Bao-Chen \cite{BC} considered recently the existence, the symmetry and the symmetry breaking region of extremal functions.
\end{Rem}

\section{Proof of Theorem \ref{thm1.1}.}
\label{section2}
Let $n\geq 2$ and $x'=(x_1,\cdot\cdot\cdot,x_{n-1},0)$ for $x=(x_i)\in\R^n$. Let $V=|x'|^{2\alpha+2}|x|^{2\beta}$ with $\alpha, \beta$ satisfying \eqref{cond1.1} with $p=2$. Consider $f=f_1 f_2$, then
\begin{align*}
\frac{{\rm div}(V\nabla f)}{f} &= \frac{{\rm div}(V\nabla f_1)}{f_1} + \frac{{\rm div}(V\nabla f_2)}{f_2} +2 V\frac{\nabla f_1\cdot\nabla f_2}{f_1f_2} =: I_1 +I_2 + I_3.
\end{align*}
Choose now $f_1=|x'|^\theta,~f_2=|x|^\lambda$, direct calculus yields that in $\R^n\setminus{\{x'=0}\}$,
\begin{align*}
\nabla V= \Big[(2\alpha+2)\frac{x'}{|x'|^2} +2\beta\frac{x}{|x|^2}\Big]V,\quad \nabla f_1= \theta |x'|^{\theta-2}x', \quad\nabla f_2=\lambda|x|^{\lambda-2}x
\end{align*}
and
\begin{align*}
\Delta f_1=\theta(n-3+\theta)|x'|^{\theta-2},\quad \Delta f_2=\lambda(n-2+\lambda)|x|^{\lambda-2}.
\end{align*}
Hence
\begin{align*}
I_1 = \frac{\nabla V\cdot \nabla f_1+ V\Delta f_1}{f_1} & = \Big[\frac{\theta(2\alpha+2)}{|x'|^2}+\frac{2\theta\beta}{|x|^2}\Big]V +V\frac{\theta(n-3+\theta)}{|x'|^2}\\
 & = \Big[ \frac{\theta(n-1+\theta+2\alpha)}{|x'|^2} +\frac{2\theta\beta}{|x|^2}\Big]V,
\end{align*}
and
\begin{align*}
I_2 =\frac{\nabla V\cdot \nabla f_2+ V\Delta f_2}{f_2}= \frac{\lambda(n+\lambda+2\alpha+2\beta)}{|x|^2}V,\quad I_3= \frac{2\theta\lambda }{|x|^2}V.
\end{align*}
Therefore we get
\begin{align}
\label{Wp=2}
\begin{split}
-\frac{{\rm div}(V\nabla f)}{f} &= -\Big[\frac{\theta(n-1+\theta+2\alpha)}{|x'|^2} + \frac{\lambda(n+2\alpha+2\beta+2\theta+\lambda)+2\beta\theta}{|x|^2}\Big] V\\
&= H_1(\theta, \lambda) \frac{V}{|x'|^2}+ H_2(\theta, \lambda) \frac{Vx_n^2}{|x'|^2|x|^2},
\end{split}
\end{align}
where
\begin{align*}
H_1(\theta, \lambda) = -\Big[\theta(n-1+\theta+2\alpha+2\beta)+ \lambda(n+2\alpha+2\beta+2\theta+\lambda)\Big]: =H(\theta)-H_2(\theta,\lambda)
\end{align*}
with
\begin{align*}
H(\theta) =-\theta(n-1+2\alpha+\theta),\quad H_2(\theta, \lambda) = \lambda(n+2\alpha+2\beta+2\theta+\lambda)+2\beta\theta.
\end{align*}

Seeing Li-Yan's estimate \eqref{ly}, we aim to find the maximum value of $H_1$ under the constraint $H_2 \geq 0$. As $\lim_{|\theta|\to \infty} H(\theta) = -\infty$, and $\lim_{|\lambda|\to \infty} H_1(\theta, \lambda) = -\infty$ uniformly for bounded $\theta$, $\max_{H_2 \geq 0}H_1$ exists.

\smallskip
It is easy to see that $\p_\theta H_1 - \p_\lambda H_1 \equiv 1$, hence $H_1$ has no critical point in $\R^2$ and $\max_{H_2 \geq 0} H_1$ is reached on the subset $\{H_2 = 0\}$, i.e. when
%~with $\lambda(n+2\alpha+2\beta+2\theta+\lambda)+2\beta\theta=0$. the last set is non empty set if and only if
\begin{align}
\label{est2.1}
\lambda(n+2\alpha+2\beta+2\theta+\lambda)+2\beta\theta=0.
\end{align}

\medskip
If $K = -4\beta(n + 2\alpha + \beta) \leq 0$, then for any $\theta \in \R$, there exists $\lambda$ such that $H_2(\theta,\lambda)=0$, because the discriminant for the quadratic equation \eqref{est2.1} of $\lambda$ satisfies
$$(n+2\alpha+2\beta+2\theta)^2 - 8\beta\theta = 4\theta^2 + 4(n+2\alpha)\theta + (n+2\alpha+2\beta)^2 \geq 0 \quad \mbox{in} \R.$$ This means that
\begin{align}\label{theta_0}
 \max_{H_2 = 0} H_1 = \max_\R H(\theta) = H(\theta_0) = \frac{(n-1+2\alpha)^2}{4}, \quad \mbox{where } \theta_0 = \frac{1-n-2\alpha}{2}.
\end{align}

\medskip
Let $K  > 0$, then \eqref{est2.1} holds true for $\lambda\neq -\beta$ and
$$
\theta=\theta(\lambda) :=-\frac{\lambda(n+2\alpha+2\beta+\lambda)}{2(\lambda+\beta)}.
$$
Clearly,
\begin{align}\label{theta_1}
\min_{\lambda<-\beta}\theta(\lambda)=\theta(\lambda_1)=-\frac{(n+2\alpha)-\sqrt K}{2}=:\theta_1, \quad \mbox{with } \lambda_1= -\beta-\frac{\sqrt K}{2}
\end{align}
and
$$\lim_{\lambda\to-\infty}\theta(\lambda)=\lim_{\lambda\to-\beta^-}\theta(\lambda)=+\infty.$$

If $K\in (0,1]$, then $\theta(\lambda_0) = \theta_0$ for $\lambda_0=-\beta-\frac{1+\sqrt {1-K}}{2}$, which yields
\begin{align}
\label{lambda0}
\max_{H_2 = 0}H_1 = H_1(\theta_0,\lambda_0)=\theta_0^2, \quad \mbox{if } K \in (0, 1].
\end{align}

Consider now $K  > 1$, we can check that
\begin{align*}
\max_{H_2 = 0, \lambda<-\beta}H_1 = H(\theta_1).
\end{align*}

By the same, for any $K > 0$, there holds
\begin{align*}
\max_{H_2 = 0, \lambda>-\beta}H_1 = H(\theta_2)\quad \mbox{with } \theta_2 = -\frac{(n+2\alpha) + \sqrt{K}}{2}.
\end{align*}
Notice that $H(\theta_2) < H(\theta_1)$ for any $K > 0$, hence for $K>1$, there holds
\begin{align}
\label{lambda1}
\begin{split}
\max_{H_2 = 0} H_1 =H_1(\theta_1,\lambda_1) = H(\theta_1) & =  -\theta_1(n-1+2\alpha+\theta_1) \\
& = \frac{n+2\alpha -\sqrt K}{2}\times\frac{n -2+2\alpha +\sqrt K}{2}\\
& = \frac{(n-1+2\alpha)^2}{4} - \frac{(\sqrt K - 1)^2}{4}.
\end{split}
\end{align}
Finally, we obtain
\begin{align*}
\max_{H_2 \geq 0} H_1 = \max_{H_2 = 0} H_1 = \frac{(n-1+2\alpha)^2-\big[\sqrt{\max(K, 1)} - 1\big]^2}{4} = C_{n, \alpha, \beta},
\end{align*}
with $K = -4\beta(n+2\alpha+\beta)$. Seeing \eqref{Wp=2} and the equality \eqref{E2},
\begin{align}
\label{p=2.1}
\int_{\mathbb{R}^n}|x'|^{2\alpha+2}|x|^{2\beta}|\nabla u|^2 dx \geq C_{n, \alpha, \beta}\int_{\mathbb{R}^n} |x'|^{2\alpha+2-2}|x|^{2\beta}u^2 dx, \quad \forall\; u \in C_c^1(\R^n \backslash\{x' = 0\}).
\end{align}

\medskip
Under the assumption \eqref{cond1.1} with $p = 2$, i.e., $2\alpha > 1-n$, $2(\alpha + \beta) > -n$, there holds (using the spherical coordinates)
\begin{align*}
\lim_{\e\to 0^+} \int_{|x'|\leq M\e, |x| \leq R} |x'|^{2\alpha}|x|^{2\beta}dx = 0, \quad \forall\; M, R > 0.
\end{align*}

For any $u \in C_c^1(\R^n)$, we consider $u_\epsilon(x) = u(x) - u(x)\eta(|x'|/\epsilon)$ for $\epsilon \in (0, 1)$, with a standard cut-off function $\eta \in C_c^1(\R)$. Applying \eqref{p=2.1} to $u_\e$ and sending $\e \to 0^+$, we can claim the estimate \eqref{p=2}.

\medskip
Now we show the sharpness of the constants $C_{n, \alpha, \beta}$ in \eqref{const1.1}. We will use the spherical coordinates for $n\geq 2$, that is
\begin{align*}
\begin{cases}
x_1 & = r\sin\varphi_1\sin\varphi_2\cdots\sin\varphi_{n-2}\sin\varphi_{n-1},\\
x_2& =r\sin\varphi_1\sin\varphi_2\cdots\sin\varphi_{n-2}\cos\varphi_{n-1},\\
&\ldots\\
%x_{n-2}& =r\sin\varphi_1\sin\varphi_2\cos\varphi_3,\\
x_{n-1}& =r\sin\varphi_1\cos\varphi_2,\\
x_n& =r\cos\varphi_1
\end{cases}
\end{align*}
where $r\in \R_+$, $\varphi_k\in [0, \pi]$ for $1\leq k\leq n-2$ if $n\geq 3$ and $\varphi_{n-1}\in [0, 2\pi]$, so
\begin{align*}
dx = r^{n-1}(\sin\varphi_1)^{n-2}(\sin\varphi_2)^{n-3}\cdot\cdot\cdot\sin\varphi_{n-2}drd\varphi_1\cdot\cdot\cdot d\varphi_{n-1}.
\end{align*}
Let $v(x)=h(s)g(r), s=|x'|$ and $r=|x|$, then
$$
\nabla v(x)=h'(s)g(r)\frac{x'}{s}+ h(s)g'(r)\frac{x}{r},
$$
and
\begin{align*}
|\nabla v(x)|^2=h'(s)^2g^2(r)+h^2(s)g'(r)^2+ 2 h(s)h'(s)g(r)g'(r)\frac{s}{r}.
\end{align*}
Hence
\begin{align*}
|x'|^{2\alpha+2}|x|^{2\beta}|\nabla v(x)|^2 & = h'(s)^2g^2(r)s^{2\alpha+2}r^{2\beta}+h^2(s)g'(r)^2 s^{2\alpha+2}r^{2\beta}+ 2 h(s)h'(s)g(r)g'(r)s^{2\alpha+3}r^{2\beta-1}\\
& = h'(s)^2g^2(r)(\sin\varphi_1)^{2\alpha+2}r^{2\beta+2\alpha+2}+h^2(s)g'(r)^2 (\sin\varphi_1)^{2\alpha+2}r^{2\beta+2\alpha+2}\\
& \quad + 2 h(s)h'(s)g(r)g'(r)(\sin\varphi_1)^{2\alpha+3}r^{2\beta+2\alpha+2}\\
& =: J_1+J_2+J_3.
\end{align*}
Denote $\Sigma=(0, \pi)$, we have
\begin{align*}
 & \quad \int_{\mathbb{R}^n} J_1 dx\\
& = \int_{\mathbb{R}_+\times\Sigma^{n-2}\times(0,2\pi)} h'(r\sin\varphi_1)^2g^2(r)r^{n+2\alpha+2\beta+1}(\sin\varphi_1)^{n+2\alpha}(\sin\varphi_2)^{n-3}\ldots\sin\varphi_{n-2}dr d\varphi_1\ldots d\varphi_{n-1}\\
%& =2\pi \int_{\mathbb{R}_+\times\Sigma^{n-2}} h'(r\sin\varphi_1)^2g^2(r)r^{n-1}(\sin\varphi_1)^{n-2}(\sin\varphi_2)^{n-3}\ldots\sin\varphi_{n-2} dr d\varphi_1\cdot\cdot\cdot d\varphi_{n-2}\\
& = \omega_{n-1}\int_{\mathbb{R}_+\times\Sigma} h'(r\sin\varphi_1)^2g^2(r)r^{n+2\alpha+2\beta+1}(\sin\varphi_1)^{n+2\alpha} dr d\varphi_1,
\end{align*}
where $\omega_{n-1}$ stands for the volume of the unit sphere in $\R^{n-1}$. For the estimates of $J_i$, we consider three subcases.

\medskip
{\bf Case $K>1$.} Seeing \eqref{lambda1}, we choose the test function $v = hg$ with $h(s)=s^{\theta_1},$ $g(r)=(r^2+\epsilon^2)^{\frac{\lambda_1}{2}} \eta(r)$ and $\eta \in C_c^1(\mathbb{R})$ a standard cut-off function. Then
\begin{align*}
h'(r\sin\varphi_1)^2g^2(r)r^{n+2\alpha+2\beta+1}(\sin\varphi_1)^{n+2\alpha} = \theta_1^2(\sin\varphi_1)^{\sqrt K-2}g^2(r)r^{2\beta-1+\sqrt K},
\end{align*}
so
\begin{align*}
\int_{\mathbb{R}^n} J_1 dx=\omega_{n-1}\theta_1^2\int_\Sigma(\sin\varphi_1)^{\sqrt K-2}d\varphi_1 \int_0^{+\infty}g^2(r)r^{2\beta-1+\sqrt K}dr.
\end{align*}
For any $\lambda > -1$, there holds
\begin{align}
\label{beta2}
\int_0^\pi (\sin s)^\lambda ds = B\Big(\frac{\lambda+1}{2},\frac{1}{2}\Big),
\end{align}
where $B(\cdot, \cdot)$ stands for Euler's Beta function
\begin{align*}
B(t,\gamma):=\int_0^1 s^{t-1}(1-s)^{\gamma-1}ds, \quad \forall\; t, \gamma > 0.
\end{align*}
We know also that
\begin{align}
\label{beta1}
B(t+1,\gamma)=\frac{t}{t+\gamma}B(t,\gamma), \quad \forall\; t, \gamma > 0.
\end{align}
% let $x=\sin^2\theta$, then we have
On the other hand, for any $\epsilon \in (0, 1)$,
\begin{align*}
\int_0^{+\infty}g^2(r)r^{2\beta-1+\sqrt K}dr& =\int_0^{+\infty} (r^2+\epsilon^2)^{\lambda_1}r^{2\beta-1+\sqrt K}\eta^2(r) dr\\
& = \int_0^{+\infty} (t^2+1)^{-\beta-\frac{\sqrt K}{2}}t^{2\beta-1+\sqrt K }\eta^2(\epsilon t) dt\\
& =-\ln\epsilon +O(1).
\end{align*}
Here and after, $O(1)$ means a quantity uniformly bounded for $\epsilon \in (0, 1)$. Indeed, we used the following fact.
\begin{lem}
\label{lem1}
Assume that $\xi \in L^1_{loc}(\R_+)$ and $\xi(s) -s^{-1} \in L^1([1, \infty))$. Let $\zeta \in C_c^1(\R)$ be a standard cut-off function, then there is $C > 0$ such that for any $\epsilon \in (0, 1)$,
\begin{align*}
\Big|\int_0^{+\infty} \xi(t)\zeta(\epsilon t) dt + \ln\epsilon\Big| \leq C.
\end{align*}
\end{lem}
There holds then
\begin{align*}
\int_{\mathbb{R}^n} J_1 dx=\omega_{n-1}\theta_1^2B\Big(\frac{\sqrt K-1}{2}, \frac{1}{2}\Big)|\ln\epsilon| +O(1), \quad \forall\; \epsilon \in (0, 1).
\end{align*}
Similarly, as
\begin{align*}
g'^2(r)=\lambda_1^2r^2(r^2+\epsilon^2)^{\lambda_1-2}\eta^2(r)+2\lambda_1r(r^2+\epsilon^2)^{\lambda_1-1}\eta\eta'(r)+(r^2+\epsilon^2)^{\lambda_1}\eta'^2(r),
\end{align*}
for $\epsilon \in (0, 1)$ we have
\begin{align*}
\int_0^{+\infty}g'^2(r)r^{2\beta+1+\sqrt K}dr& = \lambda_1^2\int_0^{+\infty} (r^2+\epsilon^2)^{\lambda_1-2}r^{2\beta+3+\sqrt K}\eta^2(r) dr\\
&\quad +2\lambda_1\int_0^{+\infty} (r^2+\epsilon^2)^{\lambda_1-1}r^{2\beta+2+\sqrt K}\eta\eta'(r) dr\\
& \quad +\int_0^{+\infty} (r^2+\epsilon^2)^{\lambda_1}r^{2\beta+1+\sqrt K}\eta'^2(r) dr\\
& = \lambda_1^2\int_0^{+\infty} (t^2+1)^{-\beta-\frac{\sqrt K}{2}-2}t^{2\beta+3+\sqrt K }\eta^2(\epsilon t) dt +O(1)\\
& = \lambda_1^2|\ln\epsilon| +O(1).
\end{align*}
We see that
\begin{align*}
\int_{\mathbb{R}^n} J_2 dx& =\omega_{n-1}\int_\Sigma (\sin\varphi_1)^{\sqrt K}d\varphi_1\int_0^{+\infty}g'(r)^2 r^{2\beta+1+\sqrt K}dr\\
%& =\omega_{n-1}B(\frac{\sqrt {2n-3}+1}{2},~\frac{1}{2})\int_0^{+\infty}g'(r)^2 r^{\sqrt {2n-3}}dr\\
& =\omega_{n-1} \lambda_1^2B\Big(\frac{\sqrt K+1}{2}, \frac{1}{2}\Big) |\ln\epsilon| +O(1),
\end{align*}
and also
\begin{align*}
\int_{\mathbb{R}^n} J_3 dx& =2\omega_{n-1}\times\theta_1\int_\Sigma (\sin\varphi_1)^{\sqrt K}d\varphi_1\int_0^{+\infty}g(r)g'(r) r^{2\beta+\sqrt K}dr\\
& = 2\omega_{n-1}\theta_1B\Big(\frac{\sqrt K+1}{2},\frac{1}{2}\Big)\int_0^{+\infty}g(r)g'(r) r^{2\beta+\sqrt K}dr\\
& = 2\omega_{n-1} \theta_1\lambda_1 B\Big(\frac{\sqrt K+1}{2},\frac{1}{2}\Big)|\ln\epsilon| +O(1).
\end{align*}
Finally, using \eqref{beta1}, we arrive at
\begin{align}
\label{est1.1}
\int_{\mathbb{R}^n} |x'|^{2\alpha+2}|x|^{2\beta}|\nabla v|^2 dx = A\omega_{n-1}B\Big(\frac{\sqrt K-1}{2},\frac{1}{2}\Big)|\ln\epsilon| +O(1),
\end{align}
where
$$A = \theta_1^2 + (\lambda_1^2+2\theta_1\lambda_1)\frac{\sqrt K-1}{\sqrt K}.$$
Recall that $\theta_1 = \frac{-(n+2\alpha) + \sqrt{K}}{2}$, $\lambda_1 = -\beta - \frac{\sqrt{K}}{2}$ and $K = -4\beta(n+2\alpha + \beta)$, so $\lambda_1^2 = 2\beta\theta_1$ and
\begin{align*}
A =  \theta_1^2 + \theta_1(2\beta+2\lambda_1)\frac{\sqrt K-1}{\sqrt K} & = \theta_1^2 + \theta_1(1-\sqrt K)\\
& = \frac{(n-1+2\alpha)^2}{4} -\frac{(\sqrt K -1)^2}{4}.
\end{align*}
Moreover,
\begin{align*}
\int_{\mathbb{R}^n}  |x'|^{2\alpha}|x|^{2\beta} v^2(x) dx & =\omega_{n-1} \int_{\mathbb{R_+}\times\Sigma }(\sin\varphi_1)^{n+2\alpha+2\theta_1-2}g^2(r)r^{n+2\alpha+2\beta+2\theta_1-1} drd\varphi_1\\&=\omega_{n-1} \int_{\mathbb{R_+}\times\Sigma }(\sin\varphi_1)^{\sqrt K-2}g^2(r)r^{2\beta+\sqrt K-1} drd\varphi_1
\\& =\omega_{n-1}B\Big(\frac{\sqrt K-1}{2}, \frac{1}{2}\Big)\int_0^{+\infty}g^2(r) r^{2\beta+\sqrt K-1} dr\\
& = \omega_{n-1} B\Big(\frac{\sqrt K-1}{2}, \frac{1}{2}\Big) |\ln\epsilon| +O(1).
\end{align*}
Therefore,
\begin{align*}
\lim_{\epsilon\to 0^+} \frac{\ds\int_{\mathbb{R}^n} |x'|^{2\alpha+2}|x|^{2\beta}|\nabla v(x)|^2 dx}{\ds\int_{\mathbb{R}^n}  |x'|^{2\alpha}|x|^{2\beta}{v^2(x)} dx}= \frac{(n-1+2\alpha)^2}{4} -\frac{(\sqrt K -1)^2}{4}.
\end{align*}

Notice that the function $h(s)=s^{\theta_1}$ is not smooth at $s= 0$, however as all involved integrals converge, we may use eventually a family of smooth functions to approximate $h$, so we omit the details. This means that for $K > 1$, the constant $\frac{(n-1+2\alpha)^2}{4} -\frac{(\sqrt K -1)^2}{4}$ is sharp to claim \eqref{p=2}.

\medskip
The analysis for other cases are similar, we will go through quickly.

\smallskip
{\bf Case $K =1$.} Here we take $v(x) = h(s)g(r)$ with $h(s) = s^{\theta_0 +\sigma}$, $g(r) = (r^2 + \e^2)^\frac{\lambda_0 -\sigma}{2}\eta(r)$, where $s = |x'|$, $r = |x|$, $\sigma > 0$ and $\lambda_0 = -\beta - \frac{1}{2}$.  Remark that \eqref{cond1.1} with $p=2$ yields
\begin{align}
\label{estK}
n + 2\alpha + \beta = \frac{n+2\alpha+2\beta}{2} + \frac{n+2\alpha}{2} > \frac{1}{2}.
\end{align}
Therefore $K =1$ implies $2\beta \in ({-1}, 0)$. We get then
\begin{align*}
\int_{\R^n} |x'|^{2\alpha}|x|^{2\beta}v^2(x) dx & = \omega_{n-1} \int_0^\pi  (\sin\varphi_1)^{2\sigma -1}d\varphi_1 \times \int_0^{+\infty} r^{2\beta + 2\sigma}(r^2 + \e^2)^{-\beta -\frac{1}{2}-\sigma}\eta^2(r)dr\\
& = \omega_{n-1} B\Big(\sigma, \frac{1}{2}\Big)|\ln\e| + O_\sigma(1).
\end{align*}
Here $\xi(r) = r^{2\beta + 2\sigma}(r^2 + \e^2)^{-\beta -\frac{1}{2}-\sigma}$ satisfies all assumptions of Lemma \ref{lem1}, and $O_\sigma(1)$ stands for a quantity uniformly bounded for $\e \in (0, 1)$ when $\sigma > 0$ is fixed. By the same, there holds
\begin{align*}
& \quad \int_{\R^n} |x'|^{2\alpha+2}|x|^{2\beta} |\nabla v|^2 dx\\
 & = \omega_{n-1} |\ln\e| \left[(\theta_0 + \sigma)^2 B\Big(\sigma, \frac{1}{2}\Big) + \big(\lambda_0^2 + 2\lambda_0\theta_0 - 2\theta_0\sigma-\sigma^2\big) B\Big(\sigma+1, \frac{1}{2}\Big)\right] + O_\sigma(1)\\
& = \omega_{n-1} |\ln\e| \left[ (\theta_0 + \sigma)^2 + \big(\lambda_0^2 + 2\lambda_0\theta_0 - 2\theta_0\sigma-\sigma^2\big)\frac{2\sigma}{2\sigma+1}\right] B\Big(\sigma, \frac{1}{2}\Big)+ O_\sigma(1).
\end{align*}
Taking first $\e\to 0^+$ and secondly $\sigma \to 0^+$, we see that $C_{n, \alpha, \beta} \leq \theta_0^2$ for $K =1$.

\medskip
{\bf Case $K <1$.} Let $v(x) = h(s)g(r)$ with $h(s) = s^{\theta_0 + \sigma}$, $g(r) = (r^2 + \e^2)^\frac{\lambda_0}{2}\eta(r)$. Remark that $\beta >-\frac{1}{2}$ by \eqref{estK}, let
$$0 < \sigma < \frac{\sqrt{1-K}}{2} \quad \mbox{and} \quad \lambda_0 = -\beta - \frac{1+\sqrt{1-K}}{2}.$$
The above choice is motivated by \eqref{lambda0}. There holds then
\begin{align*}
\int_{\R^n} |x'|^{2\alpha}|x|^{2\beta}v^2(x) dx & = \omega_{n-1} \int_0^\pi  (\sin\varphi_1)^{2\sigma -1}d\varphi_1 \times \int_0^{+\infty} r^{2\beta + 2\sigma}(r^2 + \e^2)^{-\beta -\frac{1+\sqrt{1-K}}{2}}\eta^2(r)dr\\
&=\omega_{n-1}B\Big(\sigma, \frac{1}{2}\Big)\e^{2\sigma-\sqrt{1-K}}\times \int_0^{+\infty} t^{2\beta+2\sigma}(t^2+1)^{-\beta-\frac{1+\sqrt{1-K}}{2}}\eta^2(\e t)dt\\
& = \omega_{n-1} B\Big(\sigma, \frac{1}{2}\Big)\e^{2\sigma-\sqrt{1-K}}\times  \big[A_1 + o_\sigma(1)\big].
\end{align*}
Here $o_\sigma(1)$ stands for a quantity tending to zero as $\e$ goes to $0$ for fixed $\sigma > 0$, and
\begin{align*}
A_1 =\int_0^{+\infty} t^{2\beta+2\sigma}(t^2+1)^{-\beta-\frac{1+\sqrt{1-K}}{2}}dt < \infty.
\end{align*}
Similarly, we get
\begin{align*}
\int_{\R^n} |x'|^{2\alpha+2}|x|^{2\beta} |\nabla v|^2 dx & = \omega_{n-1} \e^{2\sigma-\sqrt{1-K}} \times  \big[A_1 + o_\sigma(1)\big](\theta_0 + \sigma)^2 B\Big(\sigma, \frac{1}{2}\Big)\\
& \quad + \omega_{n-1} \e^{2\sigma-\sqrt{1-K}}\times  \big[A_2 + o_\sigma(1)\big] B\Big(\sigma+1, \frac{1}{2}\Big)\\
& = \omega_{n-1} \e^{2\sigma-\sqrt{1-K}} \Big[ A_1 \theta_0^2 + O(\sigma) + o_\sigma(1)\Big] B\Big(\sigma, \frac{1}{2}\Big)
\end{align*}
with
\begin{align*}
A_2 & =\lambda_0^2 \int_0^{+\infty} t^{2\beta+2\sigma+4}(t^2+1)^{-\beta-\frac{1+\sqrt{1-K}}{2}-2}dt\\
& \quad + 2\lambda_0(\theta_0+\sigma) \int_0^{+\infty} t^{2\beta+2\sigma+2}(t^2+1)^{-\beta-\frac{1+\sqrt{1-K}}{2}-1}dt < \infty.
\end{align*}
Taking first $\e\to 0^+$ and secondly $\sigma \to 0^+$, we see that $C_{n, \alpha, \beta} \leq \theta_0^2$ for $K <1$. \qed

%{\bf Case $K \leq 0$.} Remark that $\beta \geq 0$ by \eqref{estK}. Let $v(x) = |x'|^{\sigma+ \theta_0} (r^2 + \e^2)^\frac{\lambda}{2}\eta(r)$, with $\sigma > 0$ and $\lambda = -\beta - \frac{1}{2} -\sigma$. As in the previous case, we can check that $C_{n, \alpha, \beta} \leq \theta_0^2$. So we are done.\qed

\medskip
\begin{Rem}
Take $\alpha = \beta = -\frac{1}{2}$ and $n \geq 3$, we have
\begin{equation}\label{1.1}
\int_{\mathbb{R}^n}\frac{|x'|}{|x|}|\nabla u|^2 dx \geq \Big[\frac{n^2-6n+6}{4}+\frac{\sqrt{2n-3}}{2}\Big]\int_{\mathbb{R}^n}\frac{u^2}{|x'||x|}dx, \quad\forall\; u\in C_c^1(\mathbb{R}^n).
\end{equation}
Here the constant $\frac{n^2-6n+6}{4}+\frac{\sqrt{2n-3}}{2}$ is sharp. In particular, the best constant for \eqref{0.2} is $\frac{2\sqrt{3} -3}{4}$.
\end{Rem}

\begin{Rem}
The estimate \eqref{1.1} is trivial when $n=2$, but an anisotropic Leray type inequality exists. Let $\B^2$ denote the open unit disc in $\R^2$, there holds
\begin{equation*}
\int_{\mathbb{B}^2}\frac{|x_1|}{|x|}|\nabla u|^2 dx \geq \frac{1}{4}\int_{\mathbb{B}^2}\frac{|x_1|}{|x|^3(\ln|x|)^2}u^2dx, \quad \forall\; u \in C_c^1(\B^2).
\end{equation*}
Here we consider $f(x) = \sqrt{-\ln |x|}$ and check that
\begin{align*}
-\frac{{\rm div}(|x_1||x|^{-1}\nabla f)}{f} = \frac{|x_1|}{4|x|^3(\ln|x|)^2} \quad \mbox{in } {\mathbb{B}^2\backslash\{0\}}.
\end{align*}
\end{Rem}

\section{Proof of Theorem \ref{thmp}} Let $V=|x'|^{(\alpha+1)p}|x|^{\beta p}$ and $f=|x'|^\gamma$, then
$$
\nabla V= \Big[(\alpha+1)p\frac{x'}{|x'|^2}+\beta p\frac{x}{|x|^2}\Big]V,
$$
and $\nabla f= \gamma |x'|^{\gamma-2}x'$, $\Delta f=\gamma(n-3+\gamma)|x'|^{\gamma-2}$.
There hold also
$$
|\nabla f|^{p-2}=|\gamma|^{p-2}|x'|^{(\gamma-1)(p-2)},\quad \nabla(|\nabla f|^{p-2})=|\gamma|^{p-2}(\gamma-1)(p-2) |x'|^{\gamma p -2\gamma-p}x'.
$$
According to the Hardy inequality \eqref{HLp}, we will calculate
\begin{align*}
-\frac{{\rm div}(V|\nabla f|^{p-2} \nabla f)}{f^{p-1}}& =-\left[\frac{{\rm div}(V \nabla f)|\nabla f|^{p-2}}{f^{p-1}}+\frac{V \nabla f\cdot\nabla( |\nabla f|^{p-2})}{f^{p-1}}\right] =: - \big(K_1+ K_2 \big).
\end{align*}
%\begin{align*} I_1=\frac{{\rm div}(V \nabla f)}{f}\times \frac{|\nabla f|^{p-2}}{f^{p-2}}\end{align*}
More precisely,
\begin{align*}
\frac{{\rm div}(V \nabla f)}{f}& =\frac{\nabla V\cdot\nabla f+V\Delta f}{f} = \Big[\frac{\gamma(n-3+(\alpha+1)p+\gamma)}{|x'|^2}+\frac{\gamma\beta p}{|x|^2} \Big]V.
\end{align*}
This yields
\begin{align*}
K_1 =\frac{{\rm div}(V \nabla f)}{f}\times \frac{|\nabla f|^{p-2}}{f^{p-2}} & = V|\gamma|^{p-2} \Big[\frac{\gamma(n-3+(\alpha+1)p+\gamma)}{|x'|^p}+\frac{\gamma\beta p}{|x'|^{p-2}|x|^2} \Big].
\end{align*}
On the other hand,
\begin{align*}
K_2= V|\gamma|^{p-2}\frac{\gamma(\gamma-1)(p-2)}{|x'|^p}.
\end{align*}
Hence
\begin{align*}
W = -\frac{{\rm div}(V|\nabla f|^{p-2} \nabla f)}{f^{p-1}}& = -V|\gamma|^{p-2}\Big[\frac{\gamma(n-3+(\alpha+1)p+\gamma)+\gamma(\gamma-1)(p-2)}{|x'|^p}+ \frac{\gamma\beta p}{|x'|^{p-2}|x|^2} \Big]\\
& =V|x'|^{-p} \Big\{ {-|\gamma|^{p-2}\gamma\big[(p-1)\gamma+n-1+p\alpha\big]} - |\gamma|^{p-2}\gamma\beta p\frac{|x'|^2}{|x|^2} \Big\}\\
& =: V|x'|^{-p}\Big[ \overline H_1(\gamma)+ \overline H_2(\gamma)\frac{|x'|^2}{|x|^2}\Big].
%& =V|\gamma|^{p-2}|x'|^{-p} \Big\{ {-\gamma\big[(p-1)\gamma+n-1+p(\beta+\alpha)\big]} + \gamma\beta p\frac{x_n^2}{|x|^2} \Big\}.
\end{align*}

We consider respectively two cases according to the sign of $\beta$.

\medskip
{\bf Case $\beta\geq 0$.} Recall that $p\alpha >1-n$. Let $n-1 + p\alpha = -p\gamma_0$, then $\gamma_0 < 0$,
\begin{align*}
\max_\R \overline H_1(\gamma)= \overline H_1(\gamma_0) = |\gamma_0|^p \quad \mbox{and} \quad \overline H_2(\gamma_0) = \beta p|\gamma_0|^{p-1} \geq 0.
\end{align*}
Thanks to \eqref{Ep} or \eqref{HLp}, we have
\begin{equation*}
\||x|^{\beta}|x'|^{\alpha+1} \nabla u\|_{L^p(\mathbb{R}^n)} \geq |\gamma_0|^p\||x|^\beta|x'|^\alpha u\|_{L^p(\mathbb{R}^n)}, \quad \forall\; u\in C_c^1(\mathbb{R}^n\backslash\{x'=0\}).
\end{equation*}
Recall that $|x|^\beta|x'|^\alpha \in L^p_{loc}(\R^n)$ under the condition \eqref{cond1.1}, similarly as for the case $p=2$, we can extend the above estimate for $u\in C_c^1(\mathbb{R}^n)$ by approximation, so $C_{n,\alpha,\beta} \geq |\gamma_0|^p.$

\medskip
Now we prove the sharpness of the above estimate. Consider $v=|x'|^{\gamma}g(r)$ with $g(r)=(r^2+\epsilon^2)^{\frac{\lambda}{2}}\eta$, a cut-off function $\eta$ and
$$\gamma=-\frac{n-1+p\alpha}{p}+\sigma = \gamma_0 + \sigma, \quad \lambda=-\beta -\frac{1}{p}-\sigma, \quad \epsilon, \sigma > 0.$$
Then for $\epsilon\in(0,1)$, applying Lemma \ref{lem1} and \eqref{beta2},
\begin{align}
\label{est3.1}
\begin{split}
\int_{\mathbb{R}^n}|x'|^{p\alpha}|x|^{p\beta}|v|^p dx& =\omega_{n-1}\int_0^\pi (\sin\varphi_1)^{n-2+(\alpha+\gamma)p}d\varphi_1\int_0^{+\infty}r^{p(\beta+\sigma)}g^p(r) dr\\
& = \omega_{n-1} B\Big(\frac{p\sigma}{2},~\frac{1}{2}\Big)\times |\ln\epsilon| +O_\sigma(1)
\end{split}
\end{align}
where $O_\sigma(1)$ stands for a quantity uniformly bounded for $\e \in (0,1)$ and fixed $\sigma \in (0, 1)$. On the other hand,
\begin{align*}
|\nabla v|^2& = |\gamma|x'|^{\gamma-2}x'g(r)+|x'|^\gamma g'(r)\frac{x}{r}\mid^2\\
& =|x'|^{2\gamma-2}\left[\gamma^2g^2(r)+|x'|^2 g'^2(r)+ 2\gamma\frac{|x'|^2}{r}gg'(r)\right] =:|x'|^{2\gamma-2} G(r).
\end{align*}
In $\B^n$ the unit ball of $\R^n$, as $\eta\equiv 1$, we have
\begin{align*}
G(r)=\gamma^2(r^2+\epsilon^2)^\lambda + \lambda^2(r^2+\epsilon^2)^{\lambda-2}|x'|^2 r^2+2\gamma\lambda(r^2+\epsilon^2)^{\lambda-1}|x'|^2.
\end{align*}
Therefore
\begin{align*}
|x'|^{p(\alpha+1)}|x|^{\beta p}|\nabla v|^p & =|x'|^{p(\alpha+1)}|x|^{\beta p}|x'|^{p(\gamma-1)}G^{\frac{p}{2}}(r)\\
& =|x'|^{p(\alpha+\gamma)}|x|^{p\beta}(r^2+\epsilon^2)^{\frac{p\lambda}{2}}\left(\gamma^2+
\lambda^2\frac{|x'|^2r^2}{(r^2+\epsilon^2)^2}+2\gamma\lambda\frac{|x'|^2}{r^2+\epsilon^2}\right)^{\frac{p}{2}}.
\end{align*}
Let $\e, \sigma \in (0, 1)$, clearly
$$\Big|\lambda^2\frac{|x'|^2r^2}{(r^2+\epsilon^2)^2}+2\gamma\lambda\frac{|x'|^2}{r^2+\epsilon^2}\Big| \leq C\frac{|x'|^2}{r^2+\epsilon^2}\quad \mbox{in } \R^n.$$
By mean value theorem, as $|x'| \leq r$, there holds
$$\left|\left(\gamma^2 + \lambda^2 \frac{|x'|^2r^2}{(r^2+\epsilon^2)^2}+2\gamma\lambda \frac{|x'|^2}{r^2+\epsilon^2}\right)^{\frac{p}{2}}-|\gamma|^p\right| \leq C\frac{|x'|^2}{r^2+\epsilon^2} \quad \mbox{in } \R^n.$$
Using spherical coordinates, we see that
\begin{align}
\label{est3.2}
\begin{split}
&\quad \int_{\B^n}|x'|^{p(\alpha+1)}|x|^{\beta p}|\nabla v|^p dx\\
& =\omega_{n-1}\int_0^\pi (\sin\varphi_1)^{n-2+p(\alpha+\gamma)} d\varphi_1\\
& \quad \times \int_0^1 r^{n-1+p(\alpha+\beta+\gamma)}(r^2+\epsilon^2)^{\frac{p\lambda}{2}}\left(\gamma^2+
\lambda^2\frac{|x'|^2r^2}{(r^2+\epsilon^2)^2}+2\gamma\lambda\frac{|x'|^2}{r^2+\epsilon^2}\right)^{\frac{p}{2}} dr\\
& =: \omega_{n-1}|\gamma|^p\int_0^\pi (\sin\varphi_1)^{p\sigma-1} d\varphi_1 \int_0^1 r^{p(\beta+\sigma)}(r^2+\epsilon^2)^{\frac{p\lambda}{2}}dr + L_1\\
& =\omega_{n-1}|\gamma|^p B\Big(\frac{p\sigma}{2},\frac{1}{2}\Big)|\ln\epsilon| + O_\sigma(1)+ L_1.
\end{split}
\end{align}
Moreover,
\begin{align}
\label{est3.3}
\begin{split}
|L_1| & \leq C\int_0^\pi (\sin\varphi_1)^{p\sigma+1} d\varphi_1 \int_0^1 r^{p(\beta+\sigma)+2}(r^2+\epsilon^2)^{\frac{p\lambda}{2}-1} dr\\
& \leq C\Big[B\Big(\frac{p\sigma}{2}+1,\frac{1}{2}\Big)|\ln\epsilon|+O_\sigma(1)\Big]\\
& = C \frac{p\sigma}{p\sigma+1}B\Big(\frac{p\sigma}{2},\frac{1}{2}\Big)|\ln\epsilon|+O_\sigma(1).
\end{split}
\end{align}

\medskip
On $2\B^n\setminus \B^n$, directly calculation gives
\begin{align*}
|\nabla v|^2=|x'|^{2\gamma-2}\left[\gamma^2g^2(r)+|x'|^2 g'^2(r)+ 2\gamma\frac{|x'|^2}{r}gg'(r)\right]\leq C|x'|^{2\gamma-2}.
\end{align*}
Consequently
\begin{align*}
|x'|^{(\alpha+1)p}|x|^{\beta p}|\nabla v|^p \leq C|x'|^{1-n+p\sigma}|x|^{\beta p},
\end{align*}
and then
\begin{align}
\label{est3.4}
\int_{2\B^n\setminus \B^n}|x'|^{(\alpha+1)p}|x|^{\beta p}|\nabla v|^p dx \leq CB\Big(\frac{p\sigma}{2},\frac{1}{2}\Big).
\end{align}
Combining \eqref{est3.2}-\eqref{est3.4}, for small enough $\sigma > 0$ and $\e\in (0, 1)$, we can claim
\begin{align*}
\int_{\mathbb{R}^n}|x'|^{(\alpha+1)p}|x|^{\beta p}|\nabla v|^p dx = \omega_{n-1}|\gamma|^pB\Big(\frac{p\sigma}{2},\frac{1}{2}\Big)|\ln\epsilon|\times \Big[1+ O(\sigma)\Big]+ B\Big(\frac{p\sigma}{2},\frac{1}{2}\Big)\times O(1) + O_\sigma(1).
\end{align*}
Tend first $\epsilon\to 0^+$, secondly set $\sigma \to 0^+$, we conclude that $C_{n,\alpha,\beta} \leq |\gamma_0|^p$ seeing \eqref{est3.1}. So $C_{n,\alpha,\beta} = |\gamma_0|^p$ for $p > 1$ and $\beta \geq 0$.

\medskip
{\bf Case $\beta<0$.} Here we take still $f(x) = |x'|^\gamma$, but rewrite
\begin{align*}
W & = V|x'|^{-p}\Big[-|\gamma|^{p-2}\gamma(n-1+(\alpha+\beta) p+(p-1)\gamma) +|\gamma|^{p-2}\gamma\beta p \frac{x_n^2}{|x|^2}\Big]\\
& =: V|x'|^{-p}\Big[\widetilde H_1(\gamma)+ \widetilde H_2(\gamma) \frac{x_n^2}{|x|^2}\Big]
\end{align*}
where
\begin{align*}
\widetilde H_1(\gamma) = -|\gamma|^{p-2}\gamma(n-1+(\alpha+\beta) p+(p-1)\gamma), \quad \widetilde H_2(\gamma) = |\gamma|^{p-2}\gamma\beta p.
\end{align*}
Denote $\tilde\gamma_0 = {-\frac{n-1+(\alpha+\beta) p}{p}}$. If $p(\alpha + \beta) > 1-n$, there hold $\tilde\gamma_0 < 0$ and
\begin{align*}
\max_\R \widetilde H_1(\gamma)= \widetilde H_1(\widetilde\gamma_0)= |\widetilde\gamma_0|^p, \quad \widetilde H_2(\widetilde\gamma_0) \geq 0.
\end{align*}
Seeing \eqref{HLp} and using approximation with functions in $C_c^1(\R^n\backslash\{x'=0\})$, we claim $C_{n,\alpha, \beta} \geq |\widetilde\gamma_0|^p$.\qed

\section{Proof of Theorem \ref{thm1.3}}
%\reset
Thanks to \eqref{condLY}, we need only to prove \eqref{speCKN} for
\begin{align}
\label{4.1}
\alpha p = \beta(p-1) + \mu, \quad \mbox{hence } \; \gamma_1p = \gamma_3(p-1) + \gamma_2 -1.
\end{align}
Let $V(x) =|x'|^{p\mu}|x|^{p\gamma_2}$ and $\vec F(x)=|x'|^{\beta-\mu}|x|^{\gamma_3-\gamma_2-1}x$, direct calculation yields that in $\mathbb{R}^n\setminus{\{x'=0}\}$,
\begin{align*}
|\vec F|=|x'|^{\beta-\mu}|x|^{\gamma_3-\gamma_2}, \quad V|\vec F|^{p-2}=|x'|^{\beta(p-2)+2\mu}|x|^{\gamma_3(p-2)+2\gamma_2}
\end{align*}
and
\begin{align*}
{\rm div}(\vec F)=(n-1+\beta-\mu+\gamma_3-\gamma_2)|x'|^{\beta-\mu}|x|^{\gamma_3-\gamma_2-1}.
\end{align*}
Hence we have, with \eqref{4.1},
\begin{align*}
W={\rm div}(V|\vec F|^{p-2}\vec F) & = V|\vec F|^{p-2} {\rm div}(\vec F) + \nabla(V|\vec F|^{p-2})\cdot \vec F\\
& = [n+p(\alpha+\gamma_1)]|x'|^{\alpha p}|x|^{\gamma_1 p}.
\end{align*}

Applying \eqref{CKNp}, as $V|\vec F|^p = |x'|^{\beta p}|x|^{\gamma_3p}$, we get immediately \eqref{speCKN} for $u \in C_c^1(\R^n\backslash\{x' = 0\}$ with \eqref{4.1}. Recall that under the condition \eqref{cond1.2}, $|x'|^{\alpha}|x|^{\gamma_1}, |x'|^\beta|x|^{\gamma_3}, |x'|^\mu|x|^{\gamma_2}\in L_{loc}^p(\mathbb{R}^n)$. Similarly as for Theorem \ref{thm1.1} and \ref{thmp}, we can extend the estimate for $u\in C_c^1(\mathbb{R}^n)$ by approximation.

\medskip

 Furthermore, if $\alpha = \beta = \mu$, then the above $\vec F(x) = |x|^{\gamma_3-\gamma_2-1}x$. Then $u_0(x) = e^{-\kappa_0^{\frac{1}{p-1}}|x|^{\gamma_3-\gamma_2+1}}$ (with suitable value $\kappa_0>0$) satisfies that the residual term in \eqref{CKNp},
$${\mathcal R}(\nabla u_0, u_0\kappa_0^\frac{1}{p-1}\vec F) \equiv 0 \quad \mbox{in }\; \R^n\backslash\{0\}.$$ Hence with $\gamma_3-\gamma_2+1 > 0$ and standard approximation, we can be convinced easily that the estimate \eqref{speCKN} is sharp. \qed

\begin{Rem}
The same proof provides also a generalization of \cite[Corollary 1.2]{DFLL} to the anisotropic case. Let $p > 1$, $p\alpha > 1-n$ and $\gamma_1p = \gamma_3(p-1) + \gamma_2 - 1$, there holds
\begin{align*}
\left\||x|^{\gamma_2}|x'|^{\alpha} \nabla u\right\|_{L^p(\mathbb{R}^n)}\left\||x|^{\gamma_3}|x'|^\alpha u\right\|^{p-1}_{L^p(\mathbb{R}^n)} \geq \frac{|n+p(\alpha+\gamma_1)|}{p}\left\||x|^{\gamma_1}|x'|^\alpha u\right\|_{L^p(\mathbb{R}^n)}^p, \quad \forall\; u \in C_c^1(\R^n\backslash\{0\}).
\end{align*}
Moreover, for either $\gamma_3 - \gamma_2 +1 > 0$, $(\gamma_2+\alpha)p \geq p-n$; or $\gamma_3 - \gamma_2 +1 < 0$, $(\gamma_2+\alpha)p \leq p-n$, the above estimate is sharp.
\end{Rem}

\section{Further remarks}
Once our paper was posted in arXiv, we are aware about a very recent work of Musina-Nazarov \cite{MN}. By interests in degenerate elliptic equation with $-{\rm div}(A(x)\nabla u)$, they were motivated to establish some Hardy inequalities with anisotropic weight, see also \cite{CMN}. In particular, among many other interesting results, by Theorem 1.3 and Theorem 5.4 in \cite{MN}, Musina-Nazarov obtained the following result. Consider $x = (y, x'')\in \mathbb{R}^k \times \mathbb{R}^{n-k}$ with $1 \le k \leq n-1$, assume that
\begin{align}\label{mn0}
n \ge 2,\quad p \ge 1, \quad k+p\alpha>0,\quad p(\alpha+\beta)>-n.
\end{align}

\smallskip
\noindent
{\bf Theorem B.} {\it Let $n, p, k, \alpha, \beta$ satisfy \eqref{mn0}, there exists positive constant $C > 0$, such that
\begin{equation}\label{mn1}
\left\||x|^{\beta}|y|^{\alpha+1} \nabla u\right\|_{L^p(\mathbb{R}^n)} \geq C\left\||x|^\beta|y|^\alpha u\right\|_{L^p(\mathbb{R}^n)}, \quad \forall\; u\in C_c^1(\mathbb{R}^n).
\end{equation}
If $p = 2$, the best constant $C_{n,k,\alpha,\beta}$ to \eqref{mn1} is given by
\begin{align}\label{mn2}
C_{n,k,\alpha,\beta}=\frac{(k+2\alpha)^2-[\sqrt{\max(K, n-k)}-(n-k)^2]}{4},
\end{align}
with $K=(n+2\alpha)^2-(n+2\alpha+2\beta)^2.$}

\medskip
Clearly, Theorem \ref{thm1.1} here corresponds to the case $k = n-1$, and our approach works for general $k \le n-1$. Indeed, let $V=|y|^{2\alpha+2}|x|^{2\beta}$, $f=|y|^\theta|x|^{\lambda}$, there holds
\begin{align*}
-\frac{{\rm div}(V\nabla f)}{f}= G_1(\theta,\lambda)\frac{V}{|y|^2}+ H_2(\theta,\lambda)V\frac{|x''|^2}{|y|^2|x|^2}, \quad \mbox{ in } \R^n\backslash\{y = 0\},
\end{align*}
where $G_1(\theta,\lambda) = -\theta(k+2\alpha+\theta) - H_2(\theta,\lambda)$, and $
H_2(\theta,\lambda)=\lambda(n+2\alpha+2\beta+2\theta+\lambda)+2\beta\theta$
is the same as in the proof of Theorem \ref{thm1.1}. Proceeding very similarly as in section \ref{section2}, we can claim the best constant $C_{n,k,\alpha,\beta}$ in \eqref{mn2}, so we omit the details.

\medskip
Similar to Theorem 1.2, we get further result for general $p \ge 1$.
\begin{thm}
Let $n,p,k,\alpha,\beta$ satisfy \eqref{mn0}, and $C_{n,k,\alpha,\beta}$ be the best constant to have \eqref{mn1}, then
\begin{align*}
C_{n,k,\alpha,\beta}=\left(\frac{k+p\alpha}{p}\right)^p, \quad \text{for any }\beta\geq 0.
\end{align*}
Moreover, $C_{n,k,\alpha,\beta}\geq \left(\frac{k+p\alpha+p\beta}{p}\right)^p$ if $\beta<0$ and $k+p(\alpha+\beta)>0.$
\end{thm}

\medskip
Here we proceed as in section 3 by choosing $V=|y|^{p(\alpha+1)}|x|^{p\beta}$ and $f=|y|^{\gamma}$, which yields
\begin{align*}
W=-\frac{{\rm div}(V|\nabla f|^{p-2}\nabla f)}{f^{p-1}}= & V|y|^{-p}\Big\{{-|\gamma|^{p-2}\gamma[(p-1)\gamma+k+p\alpha]} -|\gamma|^{p-2}\gamma\beta p\frac{|y|^2}{|x|^2}\Big\}.
\end{align*}
Using $k$ instead of $n-1$, the analysis is very similar to that for Theorem \ref{thmp}, so we omit the details.

\medskip
The $L^p$ Caffarelli-Kohn-Nirenberg type inequalities of Theorem \ref{thm1.3} can also be extended to more general weights $|x|^{\gamma_i}|y|^{\alpha_i}$ with $y \in \R^k$, we leave it to the interested readers.

\bigskip
\noindent
{\bf Acknowledgements.} The authors would like to thank Professors Yanyan Li and Jingbo Dou for sending us respectively the references \cite{BC} and \cite{MN}. The authors are partially supported by NSFC (No.~12271164) and Science and Technology Commission of Shanghai Municipality (No.~22DZ2229014). The authors are also truly grateful to the anonymous referees for their thorough reading and valuable comments.

\smallskip
\noindent \textbf{Data availability.} Data sharing is not applicable to this work as no new data
were created or analyzed in the study.


\begin{thebibliography}{99}
%\bibitem{AGS} Adimurthi; M. Grossi; S. Santra, Optimal Hardy-Rellich inequalities, maximum principle and related eigenvalue problem. J. Funct. Anal. (2006) 240: 36-83.
\bibitem{BC}J. Bao and X. Chen, On the anisotropic Caffarelli-Kohn-Nirenberg type inequalities: existence, symmetry breaking region and symmetry of
  extremal functions. To appear in Commun. Contemp. Math. (2024).

\bibitem{CMN} G. Cora, R. Musina and A.I. Nazarov, Hardy type inequalities with mixed weights in cones. Annali
della Scuola Normale Superiore (2024), Doi.org:10.2422/2036-2145.202305\_013.

\bibitem{DFLL} A. Do, J. Flynn, N. Lam and G. Lu, $L^p$-Caffarelli-Kohn-Nirenberg inequalities and their stabilities. Preprint arXiv2310.07083 (2023).

\bibitem{DLL} N. Duy, N. Lam and G. Lu, $p$-Bessel pairs, Hardy's identities and inequalities and Hardy-Sobolev inequalities with monomial weights. J. Geom. Anal. 32(4) (2022), paper No. 109.

\bibitem{FLL} J. Flynn, N. Lam and G. Lu, $L^p$-Hardy identities and inequalities with respect to the distance and mean distance to the boundary. Calc. Var. Partial Differential Equations 64(1) (2025), paper No. 22.

\bibitem{GM} N. Ghoussoub and A. Moradifam, Bessel pairs and optimal Hardy and Hardy-Rellich inequalities. Math. Ann. 349(1) (2011), 1-57.

\bibitem{HY} X. Huang and D. Ye, First order Hardy inequalities revisited. Commun. Math. Res. 38(4) (2022), 535-559.

\bibitem{LY} Y. Li and X. Yan, Asymptotic stability of homogeneous solutions of incompressible stationary Navier-Stokes equations. J. Diff. Equa. 297 (2021), 226-245.

\bibitem{LY2} Y. Li and X. Yan, Anisotropic Caffarelli-Kohn-Nirenberg type inequalities. Adv. Math. 419 (2023), paper No. 108958.

\bibitem{MN} R. Musina and A.I. Nazarov, Hardy type inequalities with mixed cylindrical-spherical weights: the general case. Preprint arXiv:2411.08585v1 (2024).

\bibitem{PT} Y. Pinchover and K. Tintarev, Ground state alternative for $p$-Laplacian with potential term. Calc. Var. Partial Differential Equations. 28(2) (2007), 179-201.

\end{thebibliography}
\end{document}